%% file: ptsline.tex
\documentclass[11pt]{amsart}
\hoffset= - 0.75 in

\usepackage{amsfonts, amsmath, amssymb, amsthm}
\usepackage{latexsym, graphicx, color}
\newtheorem{thm}{Theorem}[section]

\newtheorem{conj}[thm]{Conjecture}
\newtheorem{defn}[thm]{Definition}

\newtheorem{lem}[thm]{Lemma}
\newtheorem{cor}[thm]{Corollary}
\newtheorem{prop}[thm]{Proposition}

\newcommand{\ZZ}{\mathbb{Z}}
\newcommand{\RR}{\mathbb{R}}

\newcommand{\T}{\mathcal{T}}
\newcommand{\PP}{\mathbb{P}}
\newcommand{\TP}{\mathbb{TP}}

\begin{document}

\title{The moduli space of $n$ tropically collinear points in 
$\RR^d$}
\author{Mike Develin}
\address{Mike Develin, American Institute of Mathematics, 360 Portage 
Ave., Palo Alto, CA 94306-2244, USA}
\date{\today}
\email{develin@post.harvard.edu}

\begin{abstract}

The tropical semiring $(\RR, {\rm min}, +)$ has enjoyed a recent renaissance, owing to its connections to
mathematical biology as well as optimization and algebraic geometry. In this paper, we investigate the space of
labeled $n$-point configurations lying on a tropical line in $d$-space, which is interpretable as the space of $n$-species phylogenetic
trees. This is equivalent to the space of $n\times d$ matrices of tropical rank two, a simplicial complex. We
prove that this simplicial complex is shellable for dimension $d=3$ and compute its homology in this case,
conjecturing that this complex is shellable in general. We also investigate the space of $d\times n$ matrices of
Barvinok rank two, a subcomplex directly related to optimization, giving a complete description of this
subcomplex in the case $d=3$.

\end{abstract}

\maketitle

\section{Introduction}

In ordinary linear algebra over a field, the space of $n$ labeled points on a line in $\RR^d$ is a rather simple one.
Consider this space of all $d\times n$ matrices whose columns lie on a line. By translating the points (an action 
corresponding to adding a constant to each row of the matrix), we can assume that the first point lies at the 
origin. Modulo this translation, this space of matrices then has an obvious cone point, namely the zero matrix. 
Removing this matrix yields the space of all matrices with first column zero whose remaining $n-1$ columns lie on a 
unique line. This line is chosen from the pencil $\PP^{d-1}$ of lines through the origin, and each point lying on the 
line has one parameter in $\RR$, so the space is just $\PP^{d-1}\times \RR^{n-1}$. 

When one is not working over a field, the story is more complicated. Our goal is to investigate the 
space of $n$ points on a line in tropical $d$-space, i.e. $d$-space over the tropical semiring $(\RR, \oplus, 
\odot)$, where tropical addition $\oplus$ is defined via $a\oplus b = {\rm min}(a,b)$, and tropical multiplication 
$\odot$ is defined via $a\odot b = a+b$. This space, viewed as a subset of $\RR^{d\times n}$ as in the usual case, 
contains some trivial behavior (translation and dilation); when this is removed, what remains is an interesting polytopal 
complex. We will use tropical geometry to give a simple, pictorial decomposition of this complex.

Aside from being of intrinsic interest, this space also has connections to another burgeoning area of mathematics,
the study of phylogenetic trees. It is often convenient to work in projective space rather than $\RR^d$; indeed, all
tropical geometry naturally takes place in $\TP^{d-1} := \RR^d/(1,\ldots,1)\RR$. Here, the tropical Grassmannian
${\rm Gr}\,(2,d)$ parametrizes the space of lines~\cite{SS}, which is equivalent to the space of phylogenetic trees
with $d$ leaves~\cite{BHV}, i.e. $d$ present-day species.

If we consider the case of $n$ points on a line in $d$-space, this represents the more general case of having $n$ 
species on a tree with $d$ leaves. In real life, of course, the data will not work out exactly, but by seeing how 
close the closest fit of all $d$-dimensional trees is to the data given to us by the $n$ species (namely the distance 
data between them), we can come up with a likelihood estimate that the species fit into a tree with few branches. For 
instance, suppose we are trying to resolve where a fossil species fits into the ancestry tree of three present-day 
species. If it is related, the four species should comprise a good approximation of four points on a line in 
three-space; if it is far-flung, they should not (though of course there will be a line in four-space which 
approximates them well.)

Furthermore, determining where on the tree it lies gives us an estimate of where it lies as far as the history goes: 
is it the common ancestor of the three species? Is it a direct ancestor of one of them but not the other two? Such 
questions can be resolved by fitting a tree to the data.

Therefore, understanding the space of such $n$-point configurations is important for the study of phylogenetic trees. The simplicial complex description
we obtain will dissect this space into simple chambers, making it easy to test the data against a possible tree of best fit in each chamber.

Another motivation for the study of the tropical semiring comes from the field of optimization. A matrix has
\textit{Barvinok rank} $k$ if it can be expressed as the (tropical) sum of $k$ matrices of tropical rank one but not as
the sum of $k-1$ such matrices~\cite{DSS}; a $d\times n$ matrix has tropical rank one if all of its rows are tropical
scalar multiples of each other, which is the same thing as saying that its $ij$-th entry is equal to $x_i+y_j$ for
some $x_1,\ldots,x_d,y_1,\ldots,y_n$. The name comes from work of Alexander Barvinok and his collaborators, who
showed that the traveling salesman problem can be solved in polynomial time if the Barvinok rank of a matrix is fixed~\cite{Barv}. The
smallest nontrivial case is Barvinok rank two, where there exists a polynomial-time recognition algorithm~\cite{CR}.

The space of $n$ points on a line turns out to be equivalent to the space of $n\times d$ matrices of
tropical rank two, where the \textit{tropical rank} of a matrix is defined to be the size of its largest tropically
nonsingular square minor. In ordinary linear algebra, the translated definitions of tropical rank and Barvinok rank are equivalent (both being equal to ordinary rank), but over the tropical semiring Barvinok rank can be much bigger. However, tropical rank is always
less than or equal to Barvinok rank~\cite{DSS}, meaning that the space of matrices with Barvinok rank two is a subset
of the space of matrices with tropical rank two, and indeed it is a subcomplex. We will investigate this subcomplex
along with the larger complex of matrices of tropical rank two. We conjecture that the larger complex is 
shellable, proving this for $d=3$; however, the subcomplex of Barvinok rank-two matrices is decidedly not, and even has torsion in its homology. We 
determine this homology for $d=3$ and give a description of this complex for general $d$.

\section{Preliminaries and generalities}

Throughout this paper, it will prove convenient to work in {\emph tropical projective space} $\TP^{d-1} := \RR^d/(1,\ldots,1)\RR$. Every tropically linear subspace $L$ has the property that if $x$ is in $L$, then so
is $x+k(1,\ldots,1)$ for all $k$; this is even true for any algebraic variety. For this reason, $\TP^{d-1}$ is the natural setting for tropical geometry,
since the tropical geometry of $\RR^d$ is merely a cylinder over the tropical geometry of $\TP^{d-1}$. Whenever we refer to homology in this paper, we
mean reduced homology.

The rank of a matrix, and therefore of the point configuration consisting of its columns, can be defined in many 
different ways. We review three relevant definitions given in~\cite{DSS}, which are all equivalent over a field but 
inequivalent over the tropical semiring.

\begin{defn}
The \emph{tropical rank} of a matrix $M$ is the largest $r$ such that $M$ has a tropically nonsingular $r\times r$ 
minor. A 
square matrix is tropically singular if the tropical determinant
\[
\bigoplus_{\sigma\in S_n} (\bigodot_{i\in [n]} M_{i\sigma(i)}) = {\rm min}_{\sigma\in S_n} (\sum_{i\in [n]} 
M_{i\sigma(i)})
\]
achieves the indicated minimum twice.
\end{defn}

The tropical rank of a matrix is equal to the dimension of the tropical convex hull of its columns, defined to be the set of all tropical linear combinations of
those columns, i.e. the dimension of the image of the matrix.  This image has lineality space $(1,\ldots,1)$, and its image in tropical projective space
is a polytopal complex. This definition of rank is combinatorially nice, and is connected to the study of phylogenetic trees: $n$ points form a tree
metric if and only if the negated distance matrix has tropical rank two, in which case the tropical convex hull of the columns of this matrix is a
realization of the tree. Like ordinary polytopes, tropical convex hulls satisfy a Farkas lemma and have a facet description. For more on tropical convex
hulls, see~(\cite{DS}, \cite{Jos}).

Another etymology for the tropical operations is as the image of the ordinary operations on a power series ring $K$
under the degree map sending a power series to its minimal exponent. This allows us to define linear subspaces as the
tropical vanishing sets of linear ideals, effecting the following definition.

\begin{defn} 
The \emph{Kapranov rank} of a matrix $M$ is the smallest $r$ such that there exists a linear ideal $I\subset
K[x_1,\ldots,x_n]$ of codimension $r$ such that each column of $M$ is in the tropical vanishing set $\T(I)$. A 
vector $a\in \RR^n$ is in the tropical vanishing set $\T(I)$ if for each $f\in I$, the leading term of $f$ with 
respect to the weight vector $(1,a)$ (i.e. weight 1 on $t$ and weights $a$ on the $x_i$'s) is not a monomial.
\end{defn}

The notion of Kapranov rank of course accompanies the study of linear subspaces. The definition of tropical 
vanishing set is the natural one under the conception of the tropical semiring as coming from this power series ring; 
another way to put it is that the tropical vanishing set of an ideal is the image under the degree map of the 
ordinary vanishing set of the ideal in $K^n$. For more on tropical algebraic 
geometry and the study of linear subspaces, see (\cite{RGST}, \cite{SS}).

Our final definition of rank arises naturally in the field of optimization. 

\begin{defn}
The \emph{Barvinok rank} of a matrix $M$ is the smallest $r$ such that $M$ is the tropical sum of $r$ 
tropically rank-one matrices. A matrix has tropical rank one if its columns (equivalently, rows) are tropical scalar 
multiples 
of each other.
\end{defn}

We will use the following more tractable reformulation of Barvinok rank.

\begin{prop}\cite{DSS}
The Barvinok rank of a matrix $M\subset \RR^{d\times n}$ is the smallest $r$ such that its columns lie in the 
tropical convex hull of $r$ points in $\TP^{d-1}$.
\end{prop}

Our goal is to investigate the space of $n$ labeled points all lying on a line in $\TP^{d-1}$. These configurations correspond to $d\times n$ matrices of 
Kapranov rank two. There are two reasons why this (rank two instead of rank one) is the natural definition of points on a line. First of all,
because of the projectivization, it is natural to investigate point configurations in $\TP^{d-1}$, where everything
has one smaller dimension. This is essentially akin to reverse homogenization; if we have $n$ points on an arbitrary
affine line in ordinary $\RR^d$, in order to get them into a linear subspace one needs to add a dimension. Another
natural reason to investigate rank two matrices is that, with no additive identity, there is no such thing as a
rank-zero matrix. The structure of rank-one matrices is uninteresting; with any of the three definitions of 
rank above, these matrices are given by $M_{ij} = x_i + y_j$ for some choice of $x_i$'s and $y_j$'s. This is because 
in tropical projective space, where all the geometry lives, these point configurations are just $n$ copies of the 
same point.

In the rank two case, we have the following result.

\begin{prop}\cite{DSS}\label{rk2eq}
The properties of having tropical rank two and Kapranov rank two are equivalent. A matrix enjoying these properties has Barvinok 
rank at least two.
\end{prop}

The set of matrices of Barvinok rank two therefore comprises a subset of the matrices with tropical or Kapranov rank two. We 
will investigate this set, which is relevant in the field of optimization, along with the larger set of matrices of 
tropical or Kapranov rank two, which are $n$-point configurations on a tropical line in $\TP^{d-1}$.

As for the computation of this set, one way to do it is via Gr\"{o}bner basis methods. The definition of tropical
rank above means that we can look at the Gr\"{o}bner fans for all $3\times 3$ determinants of a $d\times n$ matrix of
indeterminates; for each, we then pick out the subcomplex where the initial form is not a monomial, and intersect
these over all $3\times 3$ minors. This immediately gives a polyhedral decomposition of this space, which will
obviously have cone point equal to the zero matrix; removing this yields a polytopal complex. We can do
the same thing to look at the space of $d\times n$ matrices of tropical rank $r$. Similarly, it follows from results
in~\cite{DSS} that the space of $d\times n$ matrices of Kapranov rank $r$ is the no-monomial subcomplex of the
Gr\"{o}bner fan of the ideal generated by all $(r+1)\times(r+1)$ determinants of a $d\times n$ matrix of indeterminates.  
Since these determinants do not form a Gr\"{o}bner basis, these two are not in general the same, although
Proposition~\ref{rk2eq} implies that these complexes are the same in our case of $r=2$.

Our goal is to use geometry in order to get a more explicit, evocative, and tractable definition of this polyhedral complex. The first step is to mod
out by translations. Adding a constant to each column does not change the point configuration or the tropical rank of a matrix. Similarly, adding a
constant to a row merely translates the point configuration, which does not change its rank. A note of caution is in order: we must at each step check
that we do this in a continuous fashion, so as to preserve the structure of the space, since obviously we can distort the space by picking a bad choice
of representatives from the cosets of the translation action. We will describe a general continuous translation process in 
Section~\ref{general}.

After normalizing in this manner, we then investigate the interesting part of the space by removing the cone point given by the zero matrix to produce 
a polytopal complex, which we denote by $T_{d,n}$. In
particular, we use tropical geometry to give a smaller and more informative polyhedral decomposition than the one given by the Gr\"{o}bner fans, and
demonstrate how this can be used to compute homology of the resulting polyhedral complex. This geometric description is useful in proving that the
complex is shellable for $d=3$; we conjecture that this is true in general.

From this description, it will also be apparent that the $d\times n$ matrices of Barvinok rank two comprise a subcomplex, which we denote by $B_{d,n}$.  
This subcomplex has interesting homology, which again is most easily computable via the use of the description emanating from tropical geometry. We give
a complete description and analysis of these complexes for $d=3$; for $d=4$ and small $n$ we present computational results.

\section{The general case}\label{general}

In this section, we give a general discourse on $d\times n$ matrices of tropical rank two, giving a general geometrically inspired polytopal
decomposition of the complex consisting of the interesting part of this space. We use this polyhedral decomposition to show that the complex is pure of
dimension $d+n-4$, as is the subcomplex $B_{d,n}$ of matrices of Barvinok rank two.

First, we describe how to construct a canonical line given a set of $n$ points on a line in $\TP^{d-1}$. Let $P$ be 
the tropical convex hull of the $n$ points; this is the union of all tropical line segments between pairs of points, 
and is a tree. By a result in~\cite{DS}, tropical line segments between two points are the concatenation of ordinary 
line segments, whose slopes are all 0/1-vectors, such that the slopes taken along the path from $x$ to $y$ form a 
descending chain in the Boolean poset. At any node of this tropical line segment, it follows immediately that the 
slopes of the one or two outgoing line segments are 0/1-vectors with disjoint supports. Note that a given such line segment, since we are working in $\TP^{d-1}$ where $(1, \ldots, 1) = 0$, has two distinct 0/1 slopes, $x$ and $-x$, where both are 0/1 vectors with complementary supports. When we talk about the slope of an outgoing line segment from a point $p$ to a point $q$, we mean the one of these which is a positive multiple (in the ordinary sense) of $q-p$. 

Therefore, given any node of the tropical convex hull, the slopes of the outgoing line segments are all 0/1-vectors
with disjoint supports in $\{1,\ldots,n\}$; if two supports intersected, then taking points in those directions would
yield a tropical line segment violating the property of the previous paragraph. At each node where the union of these
supports is not $\{1,\ldots,n\}$, we add an outgoing leaf to infinity in the positive $x_i$-direction for every $i$
not in the union. (Note that if an original point is on a coordinate leaf, this has the effect of extending that
leaf.)

\begin{figure}
\includegraphics{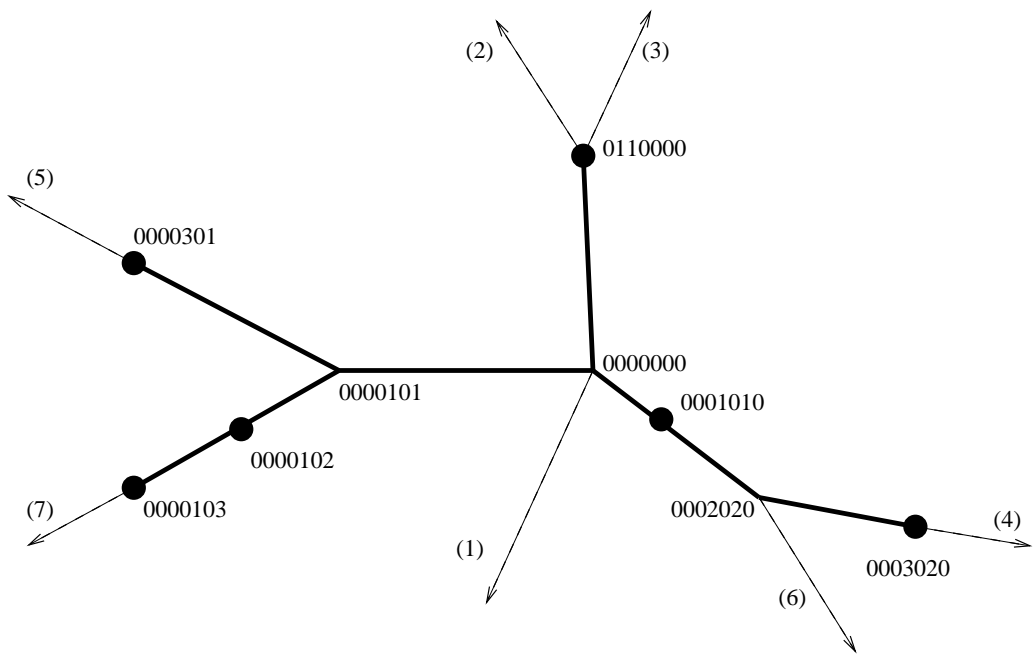}
\caption{\label{canline}
Construction of a canonical line from six points in $\TP^7$ lying on a line. The coordinate leaves are labeled with their 
direction; the points as well as the other nodes of the tree comprising their convex hull are labeled with their coordinates.}
\end{figure}

We claim that this produces a tropical line, as in Figure~\ref{canline}. This follows immediately from the fact that it produces a tree with
leaves heading to infinity in the coordinate directions, and 0/1 slope vectors which satisfy the zero tension
condition. Any such tree is in fact a tropical line~\cite{SS}. We call this the line generated by the $n$ points in question.

We now employ our translation reduction. The leaf in the direction of the first coordinate has an endpoint; we translate (the
line and) the point configuration so that this point is sent to $(0,\ldots, 0)$. This action is easily checked to be continuous,
since continuously moving the points moves this tropical line continuously. Similarly, the zero matrix consisting of all $n$
points at the origin is clearly a cone point for this set of representatives for the cosets under translation, since to find a
line containing the $k$-dilates of a point set, we need merely to dilate the line by a factor of $k$. Thus, modding out by this
translation action yields a polyhedral complex.

\begin{figure}
\includegraphics{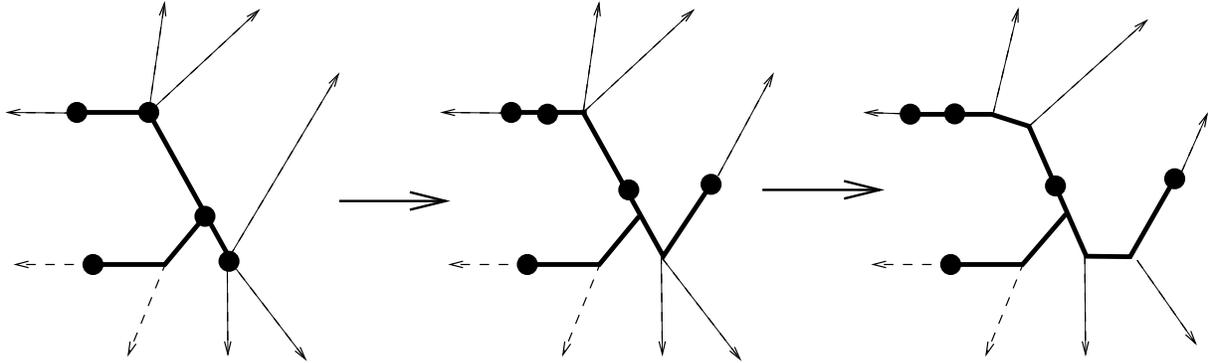}
\caption{\label{splittree}
How to find a cone of maximal dimension containing an arbitrary cone. The first step involves moving points off vertices; the 
second step involves splitting the tree at its non-trivalent nodes.}
\end{figure}

The faces of this complex will be enumerated as follows. First, pick a tree with $d$ labeled leaves (the coordinate directions.)
We fix $(0,\ldots,0)$ to be the endpoint of the leaf in the first coordinate direction. Then, for each of
the $n$ points, pick a combinatorial place on the tree for it: on a leaf, on an interior segment, or at a node. To obtain a valid face, we require that the $n$ points
must regenerate the prescribed tree; this is a combinatorial condition, stating that every internal node of the tree must lie on a segment between two of the $n$ selected points.
This setup of a fixed tree with labeled leaves and fixed combinatorial 
locations for the points will be the combinatorial object corresponding to our cone. See Figure~\ref{splittree} for some examples.

The corresponding cone consists of all $n$-point configurations which produce this particular combinatorial tree.  This cone is described by
several parameters: the lengths of the internal edges of the tree (which fixes the nodes), and for each point on a leaf or an
interior segment, its distance from an adjacent node. These cones partition the space, since each $n$-point 
configuration generates a tree via the process described above.

There is then a linear isomorphism from the set of feasible points in the parameter space to the space of all $n$-point configurations which are in
the described cone; there is clearly a linear map, and since each interior segment is actually in the convex hull of
the point configuration, changing its length will change the point configuration, so the map is injective. (This is why we had to carefully
augment the partial tree given by the tropical convex hull of the points, as opposed to taking any line containing
them.) 

Feasible points in the parameter space satisfy the following inequalities: 
for each point on an interior segment $vw$, its distance from $v$ must be 
at least zero and at most the length of $vw$, and for each point on 
a leaf, its distance from its incident vertex must be at least zero. 
This polyhedron is affinely isomorphic to the corresponding face of the 
complex. It is full-dimensional inside the parameter space, whose 
dimension is equal to the number of points on segments or leaves plus the 
number of interior edges.

We can compute the facets of this cone simply: they are the results when an inequality in the parameter space is set 
to zero. Setting the length of an interior segment to zero defines a facet only when there are no points on that 
segment; this corresponds to contracting an interior edge of the configuration. The other facets consist of sliding a
point on the interior segment $vw$ to either $v$ or $w$.

From this, we can prove easily and satisfyingly that the complex is pure. This is also implied by the classical
theorem of Bieri and Groves~\cite{BG}, which states that the tropical vanishing set $\T(I)$ is a pure polyhedral
complex if $I$ is prime (here $I$ is the ideal of all $3\times 3$ subdeterminants of a $d\times n$ matrix of
indeterminates); however, that proof is algebraic and not enlightening regarding the complex at hand. We can
construct a maximal cone which has a given cone as a face as follows. First, slide each point at a vertex onto one of
the adjacent segments; choose a leaf when the vertex is a leaf of $P$, and a segment in $P$ otherwise. This cone clearly has the original cone as a face. Let $P$ be
the convex hull of the resulting configuration.

We now claim that the only internal nodes of the tree generated by $P$ are nodes incident on two ordinary line segments from $P$. Indeed,
if a leaf of $P$ were at an internal node, it would be incident on a leaf of the tree by construction, and we would have moved it
onto that leaf in the first step. Consequently, each internal node of this tree has at least two line segments of $P$ adjacent to
it. If this node is not trivalent, insert a bridge into it to split the node into two parts in such a way that
not all the incoming line segments of $P$ are on the same side (see Figure~\ref{splittree});  then the resulting point configuration generates the split tree.
Letting this bridge go to length yields the original cone, so we have constructed a cone with the original cone as a face. Keep
doing this until the tree is trivalent; the final cone will again have the original cone as a face. This process is depicted in 
Figure~\ref{splittree}.

Since its tree is trivalent, this cone it has $d-3$ interior edges. Its dimension is therefore equal to $d-3$ plus 
the number of points on an interior segment or a leaf, which is all $n$ of them. So the dimension of the cone is 
equal to $d+n-3$, the maximum possible for any cone in our decomposition, and it has the original cone as a face. 
Therefore the polyhedral complex is pure of dimension $d+n-3$, and the polytopal complex is pure of dimension 
$d+n-4$. It is easy to check that our operations do not change the property of having Barvinok rank two (which means 
that the points are all in the tropical convex hull of two points, which can be taken to be two of the original 
points), and so by the same reasoning, the subcomplex of point configurations of 
Barvinok rank two is also pure (as a polytopal complex) of dimension $d+n-4$.

The polytopal complex we have given is much smaller than the complex given by the Gr\"{o}bner decomposition. Considering only the facets, 
the Gr\"{o}bner decomposition distinguishes the order of all the points on the interior segments, as well as the 
order of all but the two furthest-out points on each coordinate leaf. Our decomposition does neither, and hence has 
far fewer facets. 

As constructed, however, our complex has the drawback of not being a simplicial complex, which makes its homology somewhat harder to compute. The number
of facets of a cone is easy to compute: it is the number of internal segments with no points, plus the number of points on leaves, plus twice the number
of points on internal segments. By some easy dimension-counting, this cone will be a simplex whenever no internal segment has multiple points. We can
refine our decomposition to make it simplicial by specifying the order of the points on internal segments. This still yields a decomposition much
smaller than the Gr\"{o}bner one, especially when $n$ is much bigger than $d$, such as in the case where $d$ is fixed and we are investigating these
complexes for arbitrary $n$.

In the remaining two sections, we will use this small, geometric decomposition to investigate the cases of $d=3$ and 
$d=4$, giving a nice combinatorial description of both $T_{3,n}$ and $B_{3,n}$ and using this to compute their 
homology and shell the former. For $d=4$, we present experimental results and give an expanded geometric picture for 
the specific case of $n=4$ and $d=4$, using this to interpret geometrically the duality given by transposing the 
matrix. 
 
\section{The two-dimensional case}\label{twod}

In this section, we investigate the case of $n$ points on a line in the tropical projective plane $\TP^2$. This is
equivalent to considering the space of $3\times n$ matrices of tropical rank two, or the points in the tropical
determinantal variety $\T(I)$, where $I$ is the ideal generated by all $3\times 3$ minors of a $3\times n$ matrix of
indeterminates. We consider this space of matrices.

Modding out by the uninteresting parts as in Section~\ref{general} leaves us with the space of $n$ points on the 
standard line with apex $(0,0,0)$, not all on the relative interior of one branch (if $n$ points all lie on the same 
branch of a line, the process in Section~\ref{general} will generate the line where one of them is the apex.)

We now investigate this polytopal complex. Each point can either be at the apex $0$, or on branch 1, 2, or 3;  
therefore, each matrix $M$ has an associated length-$n$ string in the alphabet $\{0,1,2,3\}$, which we denote by
$\phi(M)$. Let $\Delta_X$ denote the topological closure of $\phi^{-1}(X)$, where $X$ is a string of length $n$ in
that alphabet. We then have the following proposition.

\begin{prop}
As $X$ ranges over all length-$n$ strings in $\{0,1,2,3\}$, excepting the strings all of whose entries are identical, 
the sets $\Delta_X$ are simplicial cones forming a simplicial complex. The facets of $\Delta_X$ are of the form $\Delta_Y$, 
where $Y$ ranges over all strings identical to $X$ except that one non-zero character has been changed to 0.
\end{prop}

\begin{proof}
This is a special case of the general procedure of Section~\ref{general}. Explicitly, $\Delta_X$ is defined by the following set 
of inequalities and equalities on the matrix $M$:

\begin{eqnarray*}
M_{ij} = 0\:{\rm if}\:i\neq X_j, \:{\rm and} \\
M_{ij} \ge 0\:{\rm if}\:i = X_j. \\
\end{eqnarray*}

The matrices $M$ with $\phi(M) = X$ are those for which the inequalities are all strict. 
The facets of $\Delta_X$ correspond to setting one of the inequalities to an 
equality, which corresponds to changing an $X_j$ to 0. This produces $\Delta_Y$, where $Y=X$ except that $Y_j=0$ 
while $X_j\neq 0$.

The cone has as many facets as dimensions, and hence is simplicial.
\end{proof}

This simplicial complex $T_{3,n}$ is easy to describe. The dimension of a cone is equal to the number of nonzero
entries (subtract one to get the dimension of the corresponding simplex.) The facets are the strings with no zeroes,
except that the strings $(1,\ldots,1)$, $(2,\ldots,2)$, and $(3,\ldots,3)$ are not present; there are $3^n-3$ of
these. The extreme rays (points of the complex) are the strings with one nonzero entry; there are $3n$ of these.
Except for the top dimension, the number of $k$-faces of the simplicial complex is equal to $3^{k+1}{n\choose
{k+1}}$.

\begin{thm}
The complex $T_{3,n}$ is shellable. 
\end{thm}

\begin{proof}
We exhibit an explicit shelling of $T_{3,n}$. Its facets are all strings of length $n$ in the alphabet $\{1,2,3\}$ 
excepting $(1,\ldots,1), (2,\ldots,2), (3,\ldots,3)$, 
and the intersection of two facets $\Delta_X$ and $\Delta_Y$ is equal to $\Delta_Z$, where $Z_i$ is equal to 0 unless 
$X_i=Y_i$, in which case it is equal to both of them.

We construct the ``snake ordering'' of the ternary strings of length $n$ recursively as follows. If $n=1$, we have 
$1<2<3$. If $n>1$, denote by $\tilde{X}$ the string $X$ minus its first letter. Then we define $X<Y$ if:

\begin{eqnarray*}
X_1 < Y_1, \:{\rm or} \\
X_1 = Y_1 \:{\rm is}\:{\rm odd}\:{\rm and}\:\tilde{X}<\tilde{Y},\:{\rm or} \\
X_1 = Y_1 \:{\rm is}\:{\rm even}\:{\rm and}\:\tilde{X}>\tilde{Y}.\\
\end{eqnarray*}

The proof is then complete with the following lemma.
\end{proof}

\begin{lem}
The snake ordering, with any subset of the facets $(1,\ldots,1)$, $(2,\ldots,2)$, and $(3,\ldots,3)$ removed, is 
a valid shelling order for the simplicial complex which remains. 
\end{lem}

\begin{proof}
The proof is by induction. For $n=1$, the lemma is trivial. For $n=2$, it is easily checked. Suppose $n>2$; we need 
to check that if $X<Y$, then there exists some $Z<Y$ such that $Z\cap Y$ is a facet of $Y$ and $X\cap Y\subset Z\cap 
Y$. 

If $X_1=Y_1$, then by the inductive hypothesis we are done, since the set of facets with $X_1=k$ for fixed $k$ are
listed in either the snake ordering or its reverse; these orderings are isomorphic, and the facets missing are some
subset of $(1,\ldots,1)$, $(2,\ldots,2)$, and $(3,\ldots,3)$.  Next, suppose $X_1<Y_1$. In general, we can define
$Z$ by setting $Z_i = Y_i$ for $i\neq 1$ and setting $Z_1$ to be anything less than $Y_1$.

The only case in which this fails is the case where $Y = (2,1,\ldots,1)$, $X_1 = 1$, and $(1,\ldots,1)$ is removed.  
In this case, since $X\neq (1,\ldots,1)$, there exists some $j>1$ with $X_j\neq Y_j$. Then define $Z$ via $Z_i = Y_i$
for $i\neq j$ and $Z_j = X_j$. Since $n>2$, $Z$ is not the all-2's vector, so it has not been removed. Furthermore,
$Z<Y$, since $Z_1=Y_1=2$ is even and $\tilde{Z} > \tilde{Y}$ because $\tilde{Y} = (1,\ldots,1)$ is minimal in the
snake ordering on $(n-1)$-strings. Therefore, $Z$ is as desired, completing the proof.
\end{proof}

We can now easily compute the  homology of $T_{3,n}$.

\begin{cor}
$H_{n-1}(T_{3,n}) = \ZZ^{2^n-3}$, and $H_i(T_{3,n}) = 0$ for $i\neq n-1$.
\end{cor}

\begin{proof}
Since the complex is shellable, from the Mayer-Vietoris sequence for reduced homology, it immediately has only top 
homology, which is free of some rank. To compute this rank, it suffices to compute the Euler characteristic of the 
complex. Up to a sign, this is equal to:

\[
-3 (-1)^n + \sum_{k=0}^n (-1)^k 3^{k}{n\choose k} = (-1)^n (-3 + 
\sum_{k=0}^n (-1)^{n-k} 3^k {n\choose k}) = (-1)^n 
(-3 + (3-1)^n) = \pm (2^n - 3),
\]

so $H_{n-1}(T_{3,n})$ is free of rank $2^n-3$ as desired.

% The second claim follows immediately from the fact that the complex is shellable. For the first claim, from the 
% Mayer-Vietoris sequence for reduced homology, it is easy to see that $H_{n-1}(T_{3,n})$ is free, and that its rank 
% is the number of simplices added in the shelling for which all of their facets are already covered. For the snake 
% ordering on all $3^n$ ternary strings of length $n$, this number is easily shown to be $2^n$ by a simple induction. 
% Removing $(3,\ldots,3)$, the final element in the snake ordering, reduces this number by 1, as $(3,\ldots,3)$ had 
% this property; similarly, removing $(2,\ldots,2)$ does not change this property on any other element, and since 
% $(2,\ldots,2)$ had this property, removing it reduces this number also by 1. Finally, removing $(1,\ldots,1)$ 
% removes 
% this property from $(2,1,\ldots,1)$, and does not affect it on any other element. So the number of elements with 
% this 
% property is $2^n-3$, and hence that is the rank of $H_{n-1}(T_{3,n})$ as desired.
\end{proof}

Thus, unlike the classical case, the space of $n$ points on a line gets exponentially complicated (at least topologically) as $n$ grows
large, even in two-space. It is worthwhile to note that this description of the simplicial complex is in general much smaller than the
natural Gr\"{o}bner description (although not for $n=3$, where it is a minimal simplicial refinement of the non-simplicial
Gr\"{o}bner-derived complex.) The facets of the Gr\"{o}bner-derived complex in general correspond to not only a distribution of the
points among the branches, but also an ordering of all but the two furthest points on each branch. This makes them far more numerous; for
example, for $n=7$, our complex has 2184 facets and 21 extreme rays, while the Gr\"{o}bner decomposition of the same space has 48510
facets and 378 extreme rays.

\section{The subcomplex of Barvinok rank two matrices}

In the $3\times n$ case, the subspace of all matrices of Barvinok rank two
(which we will call $B_{3,n}$) exhibits quite different behavior.  As in
the case of tropical rank two, the property of having Barvinok rank two is
translation-independent, so as before we can reduce to the case of $n$
points on the standard line. Similarly, dilation does not change the
Barvinok rank of a matrix, so we can pass to the simplicial complex
description. A matrix has Barvinok rank two if and only if the
corresponding point configuration is contained in the convex hull of two
points; in our simplicial complex description, this is equal to the union
of the simplicial cones corresponding to strings including only $0$ and
two of $\{1,2,3\}$, or the union of the facets corresponding to length-$n$
ternary strings containing only two distinct symbols.

This subcomplex naturally breaks up into three parts, one for each pair of
symbols; we will call these $C_{1,2}$, $C_{2,3}$, and $C_{1,3}$. All three
of these parts are isomorphic. $C_{1,2}$ is the union of facets
corresponding to all length-$n$ strings of 1's and 2's, except for the
all-1's and all-2's vectors; this object is, as a simplicial complex,
equal to the boundary of the $n$-dimensional crosspolytope with the
interiors of two opposite facets removed. These three crosspolytopes,
which are homotopy equivalent to $S_{n-2}\times I$, are spliced together
along the boundaries of the missing facets, so we have for instance
$C_{1,2}\cap C_{1,3}$ equal to the boundary of the would-be facet
$(1,\ldots,1)$, namely $\cup_{j=1}^n \Delta_{X^j}$, where $X^j_i =
1-\delta_{ij}$. We will denote the homological sum $\sum (-1)^j
\Delta_{X^j}$ by $[1]$, and we will similarly denote the would-be
homological boundary of the missing facets $(2,\ldots,2)$ and
$(3,\ldots,3)$ by $[2]$ and $[3]$ respectively.

Computing the homology of this object is not too hard via a Mayer-Vietoris sequence. $C_{1,2}$, $C_{1,3}$, and 
$C_{2,3}$ have homology only in dimension $n-2$, while the intersection $C_{1,3}\cap C_{1,2}$ also has 
homology only in dimension $n-2$, and the intersection $C_{2,3}\cap (C_{1,2}\cup C_{1,3})$ only has homology in 
dimensions 0 and $n-2$. Consequently, by the Mayer-Vietoris sequence, $H_1(B_{3,n})=\ZZ$, and the only other possible 
homology is in dimensions $n-2$ and $n-1$.

To compute these homologies, we need to compute the cycles in $C_{1,2}$, $C_{2,3}$, and  $C_{1,3}$; all computations 
are essentially the same. The homology of $C_{1,2}$ is generated by the cycle $[1]$. For a facet in $C_{1,2}$, which 
corresponds to a string of 1's and 2's, define $\text{sgn}(F) = (-1)^r$, where $r$ is the number of $1$'s in the 
string. Then it is easy to see that the homological boundary $\partial(\sum \text{sgn}\,(F)\, F)$ is precisely 
$[1]+(-1)^n [2]$. 

We must now consider two cases. First, suppose $n$ is even; then $[1]+[2]$ is a boundary, as is $[2]+[3]$ and 
$[1]+[3]$. This immediately implies that 
\[
H_{n-2}(B_{3,n}) = (\ZZ[1]+\ZZ[2]+\ZZ[3])/([1]+[2], [1]+[3], [2]+[3]) \cong
\ZZ/2\ZZ. 
\]

Meanwhile, $H_{n-1}$ is given by the $(n-1)$-cycles whose boundary is zero, i.e. $\{(x,y,z)\mid x([1]+[2]) + 
y([2]+[3]) + z([1]+[3]) = 0\}$. However, this is just the zero vector $(x,y,z)$, so $H_{n-1}(B_{3,n}) = 0$ for $n$ 
even.

Similarly, when $n$ is odd, the boundaries are $[1]-[2]$, $[2]-[3]$, and $[3]-[1]$. We then have 
\begin{eqnarray*}
H_{n-2}(B_{3,n}) & = &
(\ZZ[1] + \ZZ[2] + \ZZ[3])/([1]-[2], [1]-[3], [2]-[3]) \cong \ZZ, {\rm and} \\
H_{n-1}(B_{3,n}) & = & 
\{x,y,z \mid x([1]-[2]) + y([1]-[3]) + z([2]-[3]) = 0\} \cong \ZZ.
\end{eqnarray*}

We have proven the following theorem about the homology of the (interesting component of) the space of $3\times n$ 
matrices of Barvinok rank two.

\begin{thm}
The nonzero homology of the simplicial complex $B_{3,n}$ is given by:

\begin{eqnarray*}
& n\:{\rm even}: & H_{n-2} = \ZZ/2\ZZ,\,H_1 = \ZZ \\
& n\:{\rm odd}: & H_{n-1} = \ZZ,\,H_{n-2} = \ZZ,\,H_1 = \ZZ \\
\end{eqnarray*}
\end{thm}

Again, our simplicial complex description is far more informative than the
Gr\"{o}bner decomposition. For $n=7$, the subcomplex of Barvinok rank two
has 378 facets in our presentation and 27720 in the Gr\"{o}bner one.

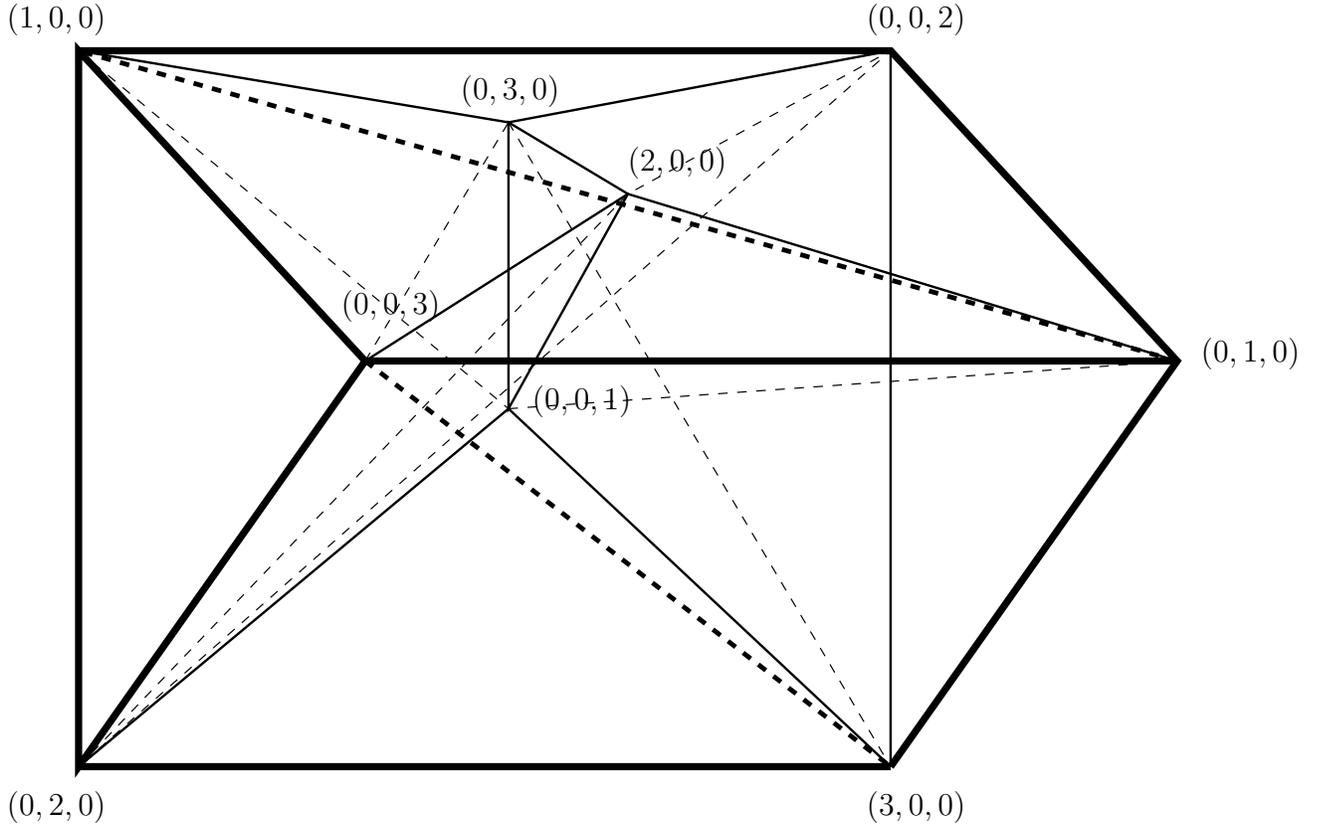
\begin{figure}
\input{schlegel.pstex_t}
\caption{\label{3x3-figure}
The polyhedral complex of $3\times 3$ matrices of tropical rank two.
}
\end{figure}

Figure~\ref{3x3-figure} shows the situation for $n=3$, that is, $3\times 3$ matrices of rank two.  The space of all $3\times 3$ matrices
of tropical rank two turns out to be the 2-skeleton of the product of two triangles, $\Delta_2\times \Delta_2$. This is represented as a
Schlegel diagram in 3-space, a triangular prism with a third parallel triangle inside it joined to the other two in the obvious way.

The Gr\"{o}bner decomposition corresponds to precisely this complex; our description corresponds to
subdividing the squares into two triangles to make the complex, the dotted lines in the figure. One can see the three-dimensional
homology $\ZZ^5$ as the five empty three-dimensional chambers. The Barvinok subcomplex consists of the union of all faces whose vertices 
only have two distinct numbers; this is the union of the nine squares in the Gr\"{o}bner decomposition, and none of the triangles. This 
union is evidently a torus, the union of three triangular prisms missing caps (or three octahedra with missing opposite facets.)

The derived cross-polytopal decomposition of $B_{3,n}$ suggests the following conjecture.

\begin{conj}
For all $m$ and $n$, the complex $B_{m,n}$ of $m\times n$ matrices of 
Barvinok rank two is a manifold.
\end{conj}

This is obviously true for $B_{3,n}$ from the decomposition into 
crosspolytopes. As we will see in the next section, $B_{4,n}$ has a more 
complicated decomposition, but the complexity of this decomposition is still independent of 
$n$, and its homology is small.

\section{The three-dimensional case}
All homology calculations in this section were done using the {\tt topaz} package of {\tt polymake}~\cite{poly}.

In this section, we discuss the case of $n$ points on a line in $\TP^3$, which corresponds to the space of all 
$4\times n$ matrices of tropical rank two. We will give a complete, nice geometric description along the lines of 
Section~\ref{twod} for the case of $4\times 4$ matrices of tropical rank two, and point out the subcomplex of 
matrices of Barvinok rank two. We will also present homology computations for small values of $n$ for these two 
complexes; as in the $n=3$ case, the Barvinok subcomplex appears to have periodic behavior with small homology, while 
the full complex has only free top homology, suggesting that the complex may be shellable.

\begin{figure}
\input{line-3space.pstex_t}
\caption{\label{line3space}
A typical line in $\TP^3$, where $s>0$ and $w,z>s$. 
}
\end{figure}
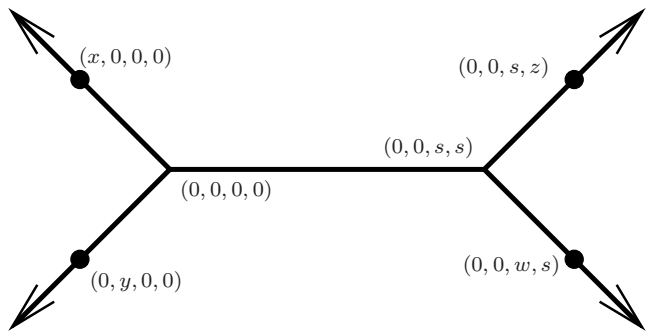

Lines in $\TP^3$ have the form shown in Figure~\ref{line3space}. Unlike in $\TP^2$, not all
lines are translates of each other. There are three distinct classes of lines, depending on which directions are
paired with each other, plus a degenerate case in which all four coordinate directions come together at a single 
point. The line shown in Figure~\ref{line3space} has directions 1 and 2 paired and directions 3 and
4 paired. In $\TP^{d-1}$, the line will look like a tree with $d$ leaves heading off in the coordinate directions;  
the slopes of the intermediate segments are as necessary for the zero tension condition. For instance, in
Figure~\ref{line3space}, at the interior point $(0,0,0,0)$, there are three segments going out, with slopes
$(1,0,0,0)$, $(0,1,0,0)$, and $(0,0,1,1)$. These slopes (taken with magnitudes so that they are 0/1-vectors) add up
to $(1,1,1,1)$, which is equal to 0 in $\TP^3$.

Deconstructing the complex $T_{4,4}$ geometrically can be done as in Section~\ref{general}, by drawing all trivalent
trees with four labeled leaves and then placing the points on the leaves and/or interior segments. There are three
such trees corresponding to the three pairings of the directions $\{1,2,3,4\}$. For each tree, we can place each of
the points $\{1,2,3,4\}$ on one of five segments. Some of these placements, however, do not yield facets. In
particular, there must be at least one point placed on one of the two leaves on each side of the tree; otherwise, the points
will not generate the tree in question. By inclusion-exclusion, this yields that the number of facets of this
decomposition is $3(5^4 - 3^4 - 3^4 + 1^4) = 1392$. The non-simplicial facets are the ones with two points on the
interior segment; there are $3(12)(4) = 144$ of these. In order to get a simplicial decomposition, we must break 
these each into two simplices by specifying the order of the interior points; this turns the 144 simplices into 288, 
for a total of 1536 simplices in our simplicial complex. There are 58 extreme rays.

\begin{figure}
\includegraphics{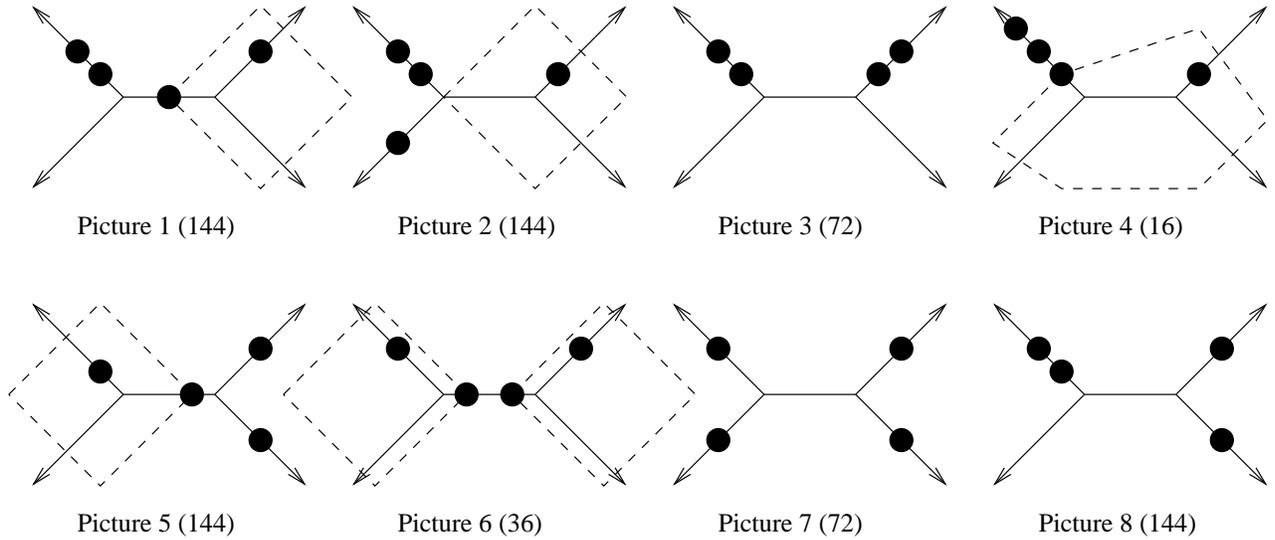}

\caption{\label{8pics4x4} The eight different ways (up to symmetry) that four points can be on a line. The boxes in Pictures 1,
2, 5, and 6 indicate that those points can be anywhere in a two-dimensional cone; the box in Picture 4 indicates that the point
inside it can be anywhere in a three-dimensional cone.}

\end{figure}

However, we can come up with an even smaller description by combining some of these simplices, again using geometry. 
The fundamental observation is the following: if two of the points are on a leaf, and one point is on the interior 
segment, then the fourth point can be anywhere in a two-dimensional orthant with apex equal to the interior point and 
directions equal to the positive coordinate directions of the two coordinates on the other side of the bridge from 
the leaf. We represent this by drawing a box around the relevant parts of the tree to denote that the point can be 
anywhere in this box. This particular situation is represented by picture 1 in Figure~\ref{8pics4x4}.

The eight pictures in Figure~\ref{8pics4x4} represent the eight different types of cones we use for a smaller
polytopal decomposition, along with the number of cones in each symmetry class (since we can assign the directions
and points arbitrarily to a given picture.) This creates a decomposition with f-vector $(34, 264, 888, 1356, 772)$
(i.e. 34 extreme rays, and 772 facets.) What we wish to highlight here is the description of the complex in terms of
valid pictures. As with our somewhat larger 1536-simplex decomposition, it is clear from the pictures what the facets
of each cone are; inside a box, they correspond to moving the point in the box to a facet of the box, while outside
the boxes they correspond to the previously described operation.

Two of our symmetry classes, pictures 4 and 6, are not simplicial. Picture 4 represents the case where three of the
points are on a coordinate leaf and the remaining point is anywhere in the three-dimensional cone given by the
closest point, which really means, from an inequality perspective, that it is inside all of the three-dimensional
cones at each of the points. We can make this simplicial by specifying which of the points is the lowest, thus
breaking this into three smaller subcones. Similarly, in Picture 6, each point in a two-dimensional cone must be
inside both two-dimensional cones with apex a point on the interior segment; distinguishing the order of the two
interior points again makes this simplicial. This introduces some extra faces of lower dimension as well, making the 
$f$-vector $(34, 264, 904, 1440, 840)$.

Computing the homology of this $T_{4,4}$ yields that $H_4(T_{4,4}) = \ZZ^{73}$ and that all other homology groups are 
zero. This, along with the shelling for $T_{3,n}$, yields the following 
conjecture.

\begin{conj}
For all $m$ and $n$, the complex $T_{m,n}$ consisting of $m\times n$ 
matrices of tropical rank two is shellable. 
\end{conj}

This conjecture is further strengthened by the computation (using similar
methods) that $T_{4,5}$ also has only top homology ($\ZZ^{301}$, to be
precise.)

\begin{figure}
\includegraphics{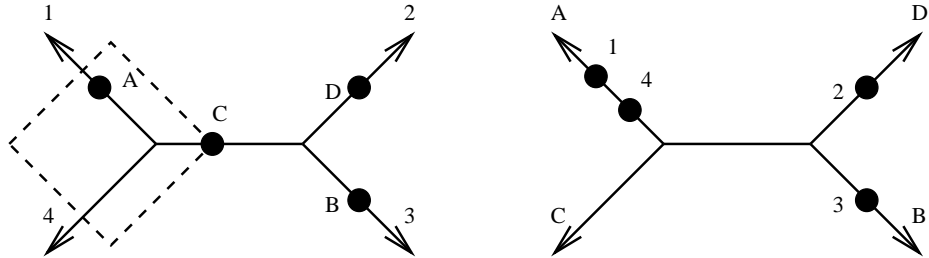}
\caption{\label{dualpic}
An example of duality. The configuration on the left, a cone in the class of Picture 5 in Figure~\ref{8pics4x4}, is dual to the 
configuration on the right, a cone in the class of Picture 8. Transposing the matrices in the left-hand cone yields the matrices 
in the right-hand cone.}
\end{figure}

Another nice feature of the geometric description given by the eight pictures in Figure~\ref{8pics4x4} is that the duality of the 
simplicial complex is clearly present. Since we are describing the space of $4\times 4$ matrices of tropical rank two, the map 
given by transposition should be an involution, and indeed we can see this in the simplicial complex. In each picture, there are 
four points and four leaves, and by looking at the relative distances between the points and the leaves, we can draw a dual 
picture where the points become the leaves and the leaves become the points. For instance, consider Figure~\ref{dualpic}. In the left-hand diagram, leaf 1 is closest to A, second closest to C, and equidistant from B and D. Leaf 2 is closest to D, then B, then C, then A. Point A (which can be anywhere in a two-dimensional cone) is closest to leaves 1 and 4, and furthest from leaves 2 and 3; the rest of the points and leaves can have their distances to members of the other category similarly ordered. In the right-hand diagram, the same orderings hold, except that the points are now represented by leaves and vice versa. (The points 1 and 4 can be on leaf A in either order in the cone represented by the right-hand diagram.) 

In this fashion, we can compute that a cone in the class of Picture 1 will be dual 
to another cone in that class, as the upper-left leaf becomes the point in the box, the two points on the upper-left leaf become 
the two leaves of the box, and so on. Similarly, the cones in Picture 2 are dual to other cones in Picture 2, the cones in 
Picture 3 are dual to the cones in Picture 6 (after those cones are split into simplices by distinguishing the order of points on 
the initial segment), the cones in Picture 4 are dual to each other (again after distinguishing the bottom point to make them 
simplicial), the cones in Picture 5 are dual to the cones in Picture 8 (as in Figure~\ref{dualpic}), and the cones in Picture 7 are 
dual to each other.

Meanwhile, the subcomplex of Barvinok rank-two configurations, $B_{4,4}$, is easily picked out; this is the union of 
the facets corresponding to pictures 1, 3, 4, and 6. Computing the homology yields some interesting results 
reminiscent of the $B_{3,n}$ case; $H_{\star}(B_{4,4}) = (0,\ZZ/2, \ZZ/2, 0, \ZZ)$. Thus it has top homology $\ZZ$, 
and torsion in dimensions 1 and 2. Computing homology for $B_{4,5}$ yields similarly striking results; we get 
$H_{\star}(B_{4,5}) = (0, \ZZ/2, 0, \ZZ, \ZZ/2, 0)$. Finally, $H_{\star}(B_{4,6}) = (0, \ZZ/2, 0, 0, \ZZ/2, 0, \ZZ)$. 

\begin{figure}
\includegraphics{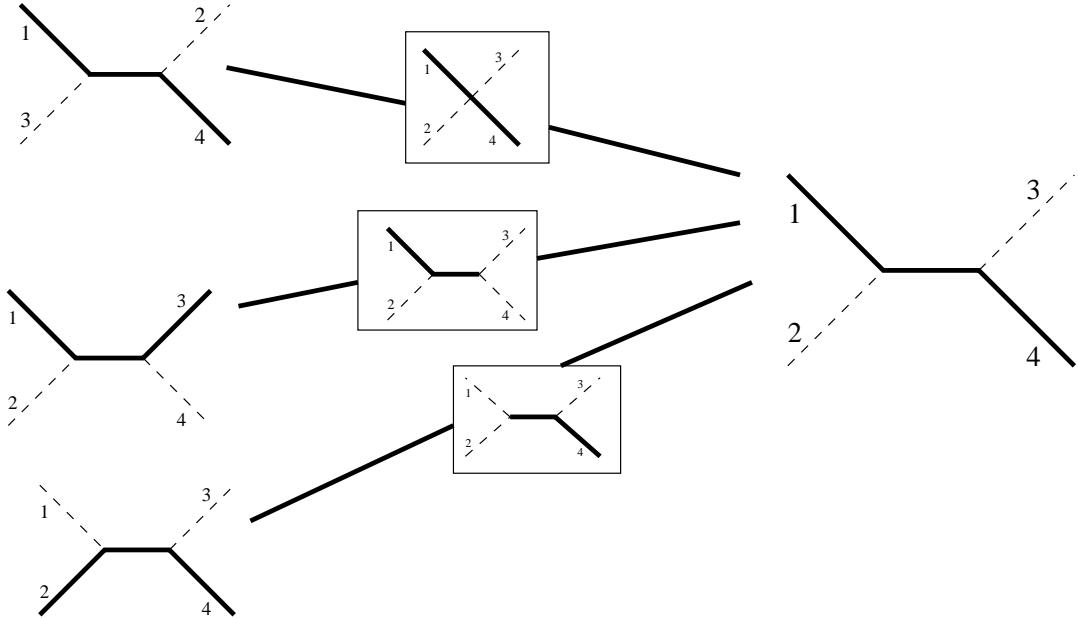}
\caption{\label{barv44}
The picture on the right-hand side represents a class of four points in $\TP^3$ of Barvinok rank two, namely the class where all
four points are in the indicated region of the indicated tree. The three pictures on the left are three other classes sharing
facets with the class on the right; the shared facet class is given by four points on the indicated region of the picture
labeling the edge.}
\end{figure}

As in the three-dimensional case, there is a suggestive way to decompose the Barvinok subcomplex $B_{4,n}$; indeed, a
similar technique works for $B_{m,n}$. The facets of $T_{4,n}$ which are in $B_{4,n}$ are those where the points are
contained in the convex hull of two points; that is, the configurations where no pair of leaves on the same side of
the bridge contain points. These decompose into twelve classes, depending on the direction of the bridge (three
choices) and a choice of a leaf from each side which is allowed to contain points (four choices.) These twelve
classes of cones intersect predictably; for instance, the class in Figure~\ref{barv44} shares facets with three other
classes as depicted. If we could now bound the homology of each of these classes and of their intersections (each of
which corresponds to a picture) independent of $n$, we could then bound the homology of the entire complex. 

Furthermore, each class has a nice combinatorial description as a complex. Its faces correspond to length-$n$ 
strings in the alphabet $\{1,2,3,A,B,C\}$, where $1$ and $2$ correspond to the leaves, $3$ corresponds to the 
interior segment, $A$ and $B$ correspond to the intersection points of the two leaves with the interior segment, and 
$C$ corresponds to the central point in the case where the interior segment has length 0. The set of strings 
corresponding to actual faces consists of all strings with the following properties:

\begin{itemize}
\item If the string contains $C$, it contains none of $\{A,B,3\}$. 
\item The string does not contain only the symbols $\{1,3,A\}$ or only the symbols $\{2,3,B\}$. 
\end{itemize}

The dimension of the cone is equal to the number of 1's, 2's, and 3's, plus one if the string does not contain any
$C$'s (i.e. corresponds to an unfused picture.) The intersection rule is also easy to give combinatorially. We simply
intersect the strings coordinate-wise, via the following intersection table, which just corresponds to finding the
largest picture which is contained in both component pictures.

$$
\begin{array}{r|r|r|r|r|r|r}
 & 1 & 2 & 3 & A & B & C \\
\hline 1 & 1 & C & A & A & C & C \\
\hline 2 & C & 2 & B & C & B & C \\
\hline 3 & A & B & 3 & A & B & C \\
\hline A & A & C & A & A & C & C \\
\hline B & C & B & B & C & B & C \\
\hline C & C & C & C & C & C & C 
\end{array}
$$

After doing the coordinate-wise intersection, if there are any C's in the resulting string, change all the A's, B's, 
and 3's to C's. The final string then corresponds to the intersection of the two cones.

This procedure works for higher $d$ as well, yielding a decomposition of
$B_{d,n}$ into simple picture classes which intersect nicely with each
other, and whose members also intersect via a similar coordinate-wise
intersection rule;  the number of such classes does not depend on $n$. In
this fashion, it seems reasonably likely that the behavior we have
observed holds in general: that the homology of $B_{d,n}$ does not
increase in complexity as $n$ gets large.

\section{Conclusion}
In this paper, we have demonstrated how the complex $T_{d,n}$ of $d\times n$ matrices with tropical rank two has a 
nice decomposition, which in general is much smaller than the Gr\"{o}bner decomposition, and which can be expressed 
in terms of pictures. In addition, the space $B_{d,n}$ of $d\times n$ matrices of Barvinok rank two is a subcomplex, 
and both complexes are pure. For dimension three, we show that $T_{d,n}$ is shellable, while $B_{d,n}$ has periodic 
and small homology; we conjecture that these are true for arbitrary $d$ and $n$. 

The above framework of pictures can also be implemented for rank higher than two. The rank-two case is the most 
important for two reasons: first of all, it has a phylogenetic interpretation, and second of all, we have a good 
description of lines. In rank higher than two, tropical and Kapranov ranks diverge, although they agree in the case 
of corank one ($n$ points on a hyperplane.) In addition, the space of $r$-planes in $\TP^{d-1}$ is in general a 
complicated one. Nonetheless, the methods we have presented here may well be applicable to such a study.

\section*{Acknowledgements}
This work was completed while the author held the AIM Postdoctoral Fellowship 2003-2008. I would also like to thank Bernd Sturmfels for 
comments on a preliminary version, and the referee for many helpful comments and revisions.

\end{document}

%% file: schlegel.pstex_t
\begin{picture}(0,0)%
\includegraphics{schlegel.pstex}%
\end{picture}%
\setlength{\unitlength}{3947sp}%
\begingroup\makeatletter\ifx\SetFigFont\undefined%
\gdef\SetFigFont#1#2#3#4#5{%
  \reset@font\fontsize{#1}{#2pt}%
  \fontfamily{#3}\fontseries{#4}\fontshape{#5}%
  \selectfont}%
\fi\endgroup%
\begin{picture}(8207,5164)(301,-4919)
\put(5701, 89){\makebox(0,0)[lb]{\smash{{\SetFigFont{12}{14.4}{\rmdefault}{\mddefault}{\updefault}{\color[rgb]{0,0,0}$(0,0,2)$}%
}}}}
\put(4201,-811){\makebox(0,0)[lb]{\smash{{\SetFigFont{12}{14.4}{\rmdefault}{\mddefault}{\updefault}{\color[rgb]{0,0,0}$(2,0,0)$}%
}}}}
\put(3151,-361){\makebox(0,0)[lb]{\smash{{\SetFigFont{12}{14.4}{\rmdefault}{\mddefault}{\updefault}{\color[rgb]{0,0,0}$(0,3,0)$}%
}}}}
\put(3601,-2311){\makebox(0,0)[lb]{\smash{{\SetFigFont{12}{14.4}{\rmdefault}{\mddefault}{\updefault}{\color[rgb]{0,0,0}$(0,0,1)$}%
}}}}
\put(5701,-4861){\makebox(0,0)[lb]{\smash{{\SetFigFont{12}{14.4}{\rmdefault}{\mddefault}{\updefault}{\color[rgb]{0,0,0}$(3,0,0)$}%
}}}}
\put(7801,-2011){\makebox(0,0)[lb]{\smash{{\SetFigFont{12}{14.4}{\rmdefault}{\mddefault}{\updefault}{\color[rgb]{0,0,0}$(0,1,0)$}%
}}}}
\put(2401,-1711){\makebox(0,0)[lb]{\smash{{\SetFigFont{12}{14.4}{\rmdefault}{\mddefault}{\updefault}{\color[rgb]{0,0,0}$(0,0,3)$}%
}}}}
\put(301,-4861){\makebox(0,0)[lb]{\smash{{\SetFigFont{12}{14.4}{\rmdefault}{\mddefault}{\updefault}{\color[rgb]{0,0,0}$(0,2,0)$}%
}}}}
\put(301, 89){\makebox(0,0)[lb]{\smash{{\SetFigFont{12}{14.4}{\rmdefault}{\mddefault}{\updefault}{\color[rgb]{0,0,0}$(1,0,0)$}%
}}}}
\end{picture}%

%% file: line-3space.pstex_t
\begin{picture}(0,0)%
\includegraphics{line-3space.pstex}%
\end{picture}%
\setlength{\unitlength}{4144sp}%
\begingroup\makeatletter\ifx\SetFigFont\undefined%
\gdef\SetFigFont#1#2#3#4#5{%
  \reset@font\fontsize{#1}{#2pt}%
  \fontfamily{#3}\fontseries{#4}\fontshape{#5}%
  \selectfont}%
\fi\endgroup%
\begin{picture}(3882,2002)(170,-1422)
\put(2448,-325){\makebox(0,0)[lb]{\smash{{\SetFigFont{8}{9.6}{\rmdefault}{\mddefault}{\updefault}{$(0,0,s,s)$}%
}}}}
\put(2926,-1028){\makebox(0,0)[lb]{\smash{{\SetFigFont{8}{9.6}{\rmdefault}{\mddefault}{\updefault}{$(0,0,w,s)$}%
}}}}
\put(2895,154){\makebox(0,0)[lb]{\smash{{\SetFigFont{8}{9.6}{\rmdefault}{\mddefault}{\updefault}{$(0,0,s,z)$}%
}}}}
\put(692,-1123){\makebox(0,0)[lb]{\smash{{\SetFigFont{8}{9.6}{\rmdefault}{\mddefault}{\updefault}{$(0,y,0,0)$}%
}}}}
\put(628,218){\makebox(0,0)[lb]{\smash{{\SetFigFont{8}{9.6}{\rmdefault}{\mddefault}{\updefault}{$(x,0,0,0)$}%
}}}}
\put(1234,-581){\makebox(0,0)[lb]{\smash{{\SetFigFont{8}{9.6}{\rmdefault}{\mddefault}{\updefault}{$(0,0,0,0)$}%
}}}}
\end{picture}%